\documentclass[10pt]{article}
\usepackage{a4}
\usepackage{epsfig}
\usepackage{amsthm,amsmath,amssymb}
\usepackage[latin3]{inputenc}

\newcommand{\n}{\noindent}

\newcommand{\ep}{\varepsilon}
\newcommand{\N}{{\rm{I\! N}}}
\newcommand{\R}{{\rm{I\! R}}}
\newcommand{\ri}{\rightarrow}
\newcommand{\riP}{\stackrel{P}{\longrightarrow}}

\newcommand{\DN}{\Delta_{j} N}
\newcommand{\DX}{\Delta_{j} X}

\newcommand{\DXd}{\Delta_{j} \tilde J_2}

\newcommand{\DXdq}{(\Delta_{i} \tilde J_2)^2}

\newcommand{\mcq}{h\ln\frac 1 h}

\newcommand{\sqrh}{\sqrt{r(h)}}
\newcommand{\rh}{r(h)}

\newcommand{\dimo}{\noindent {\it Proof.  }}

\newtheorem{Theorem}{Theorem}[section]
\newtheorem{Corollary}[Theorem]{Corollary}

\newtheorem{Lemma}[Theorem]{Lemma}
\newtheorem{Remark}[Theorem]{Remark}
\newtheorem{Proposition}[Theorem]{Proposition}
\newtheorem{example}[Theorem]{Example}
\newcommand{\ba}{\begin{array}{c}}
\newcommand{\ea}{\end{array}}
\newcommand{\bteo}{\begin{Theorem} }
\newcommand{\eteo}{\end{Theorem} }
\newcommand{\bcor}{\begin{Corollary} }
\newcommand{\ecor}{\end{Corollary} }
\newcommand{\bprop}{\begin{Proposition} }
\newcommand{\eprop}{\end{Proposition} }
\newcommand{\blem}{\begin{Lemma} }
\newcommand{\elem}{\end{Lemma} }
\newcommand{\brem}{\begin{Remark} }
\newcommand{\erem}{\end{Remark} }
\newcommand{\beqlab}{\begin{equation} \label }
\newcommand{\bex}{\begin{example}}
\newcommand{\eex}{\end{example}}
\newcommand{\beq}{\begin{equation} }
\newcommand{\eeq}{\end{equation} }
\newcommand{\bitem}{\begin{itemize}}
\newcommand{\eitem}{\end{itemize}}
\newcommand{\benum}{\begin{enumerate}}
\newcommand{\eenum}{\end{enumerate}}

\begin{document}
\title{DIFFUSION COVARIATION AND CO-JUMPS IN BIDIMENSIONAL ASSET PRICE PROCESSES
 WITH STOCHASTIC VOLATILITY AND INFINITE ACTIVITY L\'EVY JUMPS
}
\author{Fabio Gobbi\footnote{Dipartimento di Matematica per le Decisioni,
Universit\`a degli Studi di Firenze}\quad and\quad Cecilia
Mancini\footnote{Dipartimento di Matematica per le Decisioni,
Universit\`a degli Studi di Firenze}}

\vspace{-2cm} \maketitle

\begin{abstract}
In this paper we consider two processes driven by diffusions and
jumps. The jump components are Lévy processes and they can both
have finite activity and infinite activity. Given discrete
observations we estimate the covariation between the two diffusion
parts and the co-jumps.  The detection of the co-jumps allows to
gain insight in
the dependence structure of the jump components and has important applications in finance.\\
Our estimators are based on a threshold principle allowing to
isolate the jumps. 
This work follows Gobbi and Mancini (2006) where the asymptotic
normality for the estimator of the covariation, with convergence
speed $\sqrt h$, was obtained when the jump components have finite
activity. Here we show that the speed  is $\sqrt h$ only when the
activity of the jump components is moderate.\footnote{ {\bf
Corresponding author} Fabio Gobbi, Dipartimento di Matematica per
le Decisioni, Universit\`a di Firenze, via C. Lombroso 6/17, 50134
Firenze,\\ fgobbi@ds.unifi.it, phone number +39
(0)55 4796809, fax number +39 (0)55 4796800}\\

\end{abstract}

\section{Introduction}
We consider two state variables evolving as follows
$$
dX^{(1)}_{t}=
a^{(1)}_{t}dt+\sigma^{(1)}_{t}dW^{(1)}_{t}+dJ_{t}^{(1)},
$$
$$
dX^{(2)}_{t}=
a^{(2)}_{t}dt+\sigma^{(2)}_{t}dW^{(2)}_{t}+dJ_{t}^{(2)},
$$
for $t\in [0,T]$, $T$ fixed, where $W^{(2)}_{t}=\rho_{t}
W^{(1)}_{t}+\sqrt{1-\rho_{t}^{2}}W^{(3)}_{t}$;
$W^{(1)}=(W^{(1)}_{t})_{t\in [0,T]}$ and $W^{(3)}=(W^{(3)}_{t})_{
t\in [0,T]}$ are independent Wiener processes. $J^{(1)}$ and
$J^{(2)}$ are possibly correlated pure jump processes. We are
interested in the separate identification of the dependence
elements of the processes $X^{(q)}$, i.e. both of the covariation
$\int_0^T\rho_{t}\sigma^{(1)}_t\sigma^{(2)}_t dt$ between the two
diffusion parts and of the {\it co-jumps}  $\Delta
J^{(1)}_{t}\Delta
J^{(2)}_{t}$, the simultaneous jumps of $X^{(1)}$ and $X^{(2)}$.\\
Given discrete equally spaced observations $X^{(1)}_{t_{j}},
X^{(2)}_{t_{j}}, \ j= 1..n,$ in the interval $[0,T]$ (with
$t_{j}=j\frac{T}{n}$), a commonly used approach to estimate
$\int_0^T\rho_{t}\sigma^{(1)}_t\sigma^{(2)}_t dt$ is to take the
sum of cross products
$\sum_{j=1}^{n}(X^{(1)}_{t_{j}}-X^{(1)}_{t_{j-1}})
(X^{(2)}_{t_{j}}-X^{(2)}_{t_{j-1}})$; however, this estimate can
be highly biased when the processes $X^{(q)}$ contain jumps; in
fact, such a sum approaches the global quadratic covariation
$[X^{(1)},X^{(2)}]_{T}=\int_{0}^{T}\rho_{t}
\sigma_{t}^{(1)}\sigma_{t}^{(2)}dt+\sum_{0\leq t\leq T}\Delta
J^{(1)}_{t}\Delta J^{(2)}_{t}$ containing also the co-jumps. It is
crucial to single out the time intervals where the jumps have not
occurred. Our estimator is based on a threshold criterion
(\cite{man05}) allowing to isolate the jump part. In particular,
we asymptotically identify when jumps larger than a given {\it
thershold} occurred in a given time interval $]t_{j-1}, t_j]$,
depending on whether the increment $|X_{t_{j}}-X_{t_{j-1}}|$ is
too big with respect to the  threshold. In Gobbi and Mancini
(2006) we derived an asymptotically unbiased estimator of the
continuous part of the covariation process as well as of the
co-jumps. More precisely, the following threshold estimator
$$
\tilde{v}^{(n)}_{1,1}(X^{(1)},X^{(2)})_{T}=\sum_{j=1}^{n}
\Delta_{j}X^{(1)}1_{\{(\Delta_{j}X^{(1)})^{2}\leq
r(h)\}}\Delta_{j}X^{(2)}1_{\{(\Delta_{j}X^{(2)})^{2}\leq r(h)\}},
$$
is a truncated version of the realized quadratic covariation and
it is shown to be consistent to $\int_{0}^{T}\rho_{t}
\sigma_{t}^{(1)}\sigma_{t}^{(2)}dt$, as the number $n$ of
observations tends to infinity. Moreover, in the case where each
$J^{(q)}$ is a {\it finite activity} jump process (i.e. only a
finite number of jumps can occur, along each path, in each finite
time interval) we show that our estimator is asymptotically
Gaussian and converges with speed $\sqrt{h}$. Here we find the
speed of convergence of the estimator of the covariation even in
the case of infinite activity jumps, which turns out to be $\sqrt
h$ only for moderate activity of the jump
processes.\\
For the literature on non parametric inference for stochastic
processes driven by diffusions plus jumps, see Gobbi and Mancini
(2006).

Applications of the theory we present here is of strong interest
in finance, in particular in financial econometrics (see e.g.
\cite{AndBolDie05}), in the framework of portfolio risk
(\cite{ConTan04}) and for hedge funds management.

An outline of the paper is as follows. In section 2 we illustrate
the framework; in section 3 we present some preliminary results in
the case where each component $J^{(q)}$ of $X^{(q)}$ has finite
activity of jump. In section 4 we deal with the more complex case
where each $J^{(q)}$ can have an {\it infinite activity} jump
component $\tilde J_2^{(q)}$ (which makes an infinite number of
jumps in each finite time interval). We assume that such component
$\tilde J_2^{(q)}$ is a L\'evy process and we show that our
estimator is consistent and we develop some preliminaries for the
asymptotic normality in the case where $\tilde J_2^{(q)}$ have
stable-like laws and the joint law is characterized by a Copula
ranging in a given class.

\section{The framework}

Given a filtered probability space
$(\Omega,\mathcal{F},(\mathcal{F}_{t})_{t\in [0,T]},P)$, let
$X^{(1)}=(X_{t}^{(2)})_{t \in [0,T]}$ and $X^{(2)}=(X_{t}^{(2)})_{
t \in [0,T]}$ be two real processes defined by
\begin{equation}\label{modello}
  \begin{array}{c}
  X_{t}^{(1)} = \int_{0}^{t}a_{s}^{(1)}ds+\int_{0}^{t}\sigma_{s}^{(1)}
  dW_{s}^{(1)}+J_{t}^{(1)},\quad t\in [0,T], \\
  \\
  X_{t}^{(2)} =
  \int_{0}^{t}a_{s}^{(2)}ds+\int_{0}^{t}\sigma_{s}^{(2)}dW_{s}^{(2)}+
  J_{t}^{(2)},\quad t\in [0,T],
\end{array}
\end{equation}
where

    \hspace{1cm}
     \begin{figure}[h]\center{
     \begin{minipage}[c]{0.9\linewidth}
     \textbf{A1}. $W^{(1)}=(W^{(1)}_{t})_{t \in [0,T]}$ and  $W^{(2)}=(W^{(2)}_{t})_{t \in [0,T]}$
    are two correlated Wiener processes,
    with
    $\rho_{t}=Corr(W^{(1)}_{t},W^{(2)}_{t})$, $t\in [0,T]$; we can
    write
    $$
    W^{(2)}_{t}=\rho_{t} W^{(1)}_{t}+\sqrt{1-\rho_{t}^{2}}\ W^{(3)}_{t},
    $$
    where $W^{(1)}$ and $W^{(3)}$ are independent Wiener processes.\\
     \end{minipage}}\end{figure}

     \begin{figure}[h] \center{\begin{minipage}[c]{0.9\linewidth}
     \textbf{A2}. The diffusion stochastic coefficients $\sigma^{(q)}=(\sigma^{(q)}_{t})_{t\in
    [0,T]}$, $a^{(q)}=(a^{(q)}_{t})_{t\in [0,T]}$, $q=1,2$, and
    $\rho=(\rho_{t})_{t\in [0,T]}$ are  adapted càdlàg.\\
    \end{minipage}}\end{figure}
\newpage
     \begin{figure}[h] \center{\begin{minipage}[c]{0.9\linewidth}
     \textbf{A3}. For $q=1,2$
$$J^{(q)}=J_{1}^{(q)}+\tilde{J}_{2}^{(q)},$$
where $J_1^{(q)}$ are finite activity jump processes
$$
J^{(q)}_{1
t}=\int_{0}^{t}\gamma_{s}^{(q)}dN_{s}^{(q)}=\sum_{k=1}^{N_{t}^{(q)}}\gamma_{\tau^{(q)}_{k}},\quad
q=1,2,
$$
where $N^{(q)}=(N_{t}^{(q)})_{t\in [0,T]}$ are counting processes
with $E[N^{(q)}_T]<\infty$; $\{\tau^{(q)}_{k},\
k=1,...,N_{T}^{(q)}\}$ denote the instants of jump of $J_1^{(q)}$
and $\gamma_{\tau^{(q)}_{k}}$ denote the sizes of the jumps
occurred at $\tau^{(q)}_{k}$. We assume
\beq\label{CondizSulleJumpSizesDiJ1}
P(\gamma_{\tau^{(q)}_{k}}=0)=0,\quad \forall \
k=1,...,N_{T}^{(q)},\ q=1,2. \eeq Denote, for each $q=1,2$,
$\underline{\gamma}^{(q)}= \min_{k=1,...,N_{T}^{(q)}}
|\gamma_{\tau^{(q)}_{k}}|$. By
condition (\ref{CondizSulleJumpSizesDiJ1}), a.s. we have $\underline{\gamma}^{(q)}>0$.\\
 \end{minipage}}

\center{\begin{minipage}[c]{0.9\linewidth}
 \textbf{A4}.
 $\tilde J_2^{(q)}$ are infinite activity L\'evy pure jump processes of small jumps,
\beq\label{FormaEsplicTildeJ2}
\tilde{J}_{2t}^{(q)}=\int_{0}^{t}\int_{|x|\leq 1}x\
\tilde{\mu}^{(q)}(dx, ds), \eeq where $\mu^{(q)}$ is the Poisson
random measure of the jumps of $\tilde J_2^{(q)}$,
$\tilde{\mu}^{(q)}(dx, ds)=\mu^{(q)}(dx, ds)- \nu^{(q)}(dx)ds$ is
its compensated measure, where $\nu^{(q)}$ is the L\'evy measure
of $\tilde J_2^{(q)}$ (see \cite{ConTan04}).
 \end{minipage}}\end{figure}

 Each $\nu^{(q)}$ has the property that
$\nu^{(q)}(\mathbb{R}-\{0\})=\infty$, which characterizes the fact
that the path of  $\tilde{J}_{2}^{(q)}$ jumps infinitely many
times on each compact time interval. $\tilde{J}_{2}^{(q)}$ is a
compensated sum of  jumps, each of which  is  bounded in absolute
value by 1, so that substantially $J^{(q)}_{1}$  accounts for the
"big" (bigger in absolute value than $\underline{\gamma}^{(q)}$)
and rare jumps  of $X^{(q)}$, while $\tilde J_2^{(q)}$ accounts
for the very frequent and small jumps.\\

\begin{Remark}\label{RemCasoLevyJumps} {\rm
If $J^{(q)}$ is a pure jump L\'evy process, it is always possible
to decompose it as
 $$
J^{(q)}=J_{1}^{(q)}+\tilde{J}_{2}^{(q)},
$$
(see \cite{ConTan04}) where $J_1$ is a compound Poisson process
accounting for the jumps bigger in absolute value than 1, $J_1$
satisfies assumption {\bf A3} and $\tilde J_2$ is as in
(\ref{FormaEsplicTildeJ2}).
}\end{Remark}

\n {\bf Notation.} $c$ denotes any constant.

\begin{figure}[h]
\hspace{1cm}
 \center{\begin{minipage}[c]{0.9\linewidth} \textbf{A5}. Let $\alpha_q$ be the  Blumenthal Getoor index of each $J^{(q)}$,
 $q=1,2$
 (see \cite{ConTan04}). Let each $\nu^{(q)}$ satisfy:
$$ \begin{array}{lc}
{\bf A5.1} & \int_{|x|\leq \ep} x^2 \nu^{(q)} (dx)
=O( \ep^{2-\alpha_q})\\
\\
{\bf A5.2} & \int_{\ep<|x|\leq 1} |x| \nu^{(q)} (dx)=O(c - c
\ep^{1-\alpha_q}).
 \ea$$
\end{minipage}}\end{figure}

Assumption \textbf{A5} is satisfied if for instance each
     $\nu^{(q)}$ has a density $f^{(q)}(x)$  behaving as $\frac{K^{(q)}(|x|)}{|x|^{1+\alpha_q}}$
     when  $x\ri 0$, where $K^{(q)}$ is a real function with $\lim\limits_{x\ri 0} K^{(q)}(x) \in \R-\{0\}$, and
     $\alpha_q$ is the Blumenthal-Getoor index of $J^{(q)}$.\\
     In particular \textbf{A5} is true for anyone of the commonly used
models (e.g. NIG, VG,
CGMY, $\alpha$-stable, GHL).\\

\n Let, for each $n$,
$\pi_{n}^{[0,T]}=\{0=t_{0,n}<t_{1.n}<\cdot\cdot\cdot<t_{n,n}=T\}$
be a partition of $[0,T]$. We assume equally spaced subdivisions,
i.e. $h_{n}:=t_{j,n}-t_{j-1,n}=\frac{T}{n}$ for every
$n=1,2,....$. Hence $h_{n}\rightarrow 0$ as $n \rightarrow
\infty$. Let $\Delta_{j,n}X$ be the increment
$X_{t_{j,n}}-X_{t_{j-1,n}}$. To
simplify notations we write $h$ in place of $h_{n}$  and $\Delta_{j}X$ in place of $\Delta_{j,n}X$.\\

\begin{figure}[h] \center{\begin{minipage}[c]{0.9\linewidth}  \textbf{A6}. We choose a deterministic function, $h\mapsto
r(h)$, satisfying the following properties
$$
\lim_{h\rightarrow 0}
    r(h)=0, \quad \lim_{h\rightarrow
    0}\frac{hlog\frac{1}{h}}{r(h)}=0.
$$
\end{minipage}}
\end{figure}

\n We denote $r(h)$ by $r_{h}$. Denote also, for each $q=1,2$,
$$ D^{(q)}_{t} = \int_{0}^{t}a^{(q)}_{s}ds+\int_{0}^{t}\sigma^{(q)}_{s}dW_{s}^{(q)},  $$
the diffusion part of $X^{(q)}$, and
$$ Y^{(q)}_{t} =  D^{(q)}_{t} + J_{1 t}^{(q)}.$$

\section{Preliminary results}

\n By the Paul L\'evy law  of the modulus of continuity of the
Brownian motion paths (see [14]), we know that the increments of
the diffusion part of each $\Delta_{j}X^{(q)}$ tend to zero at
speed $\sqrt{h \ln \frac 1 h}$. This is the key point to
understand when an increment $\Delta_{j}X^{(q)}$ is likely to
contain some jumps. In fact if, for small $h$, $|\DX^{(q)}|>r_h >
\sqrt{h \ln \frac 1 h}$, then or some jumps of $J_1^{(q)}$
occurred, or some jumps of $\tilde J_2^{(q)}$ larger than
$2\sqrt{r_h}$ occurred (Mancini, 2005). In Gobbi and Mancini
(2006) we obtain the following consequences.

\begin{Remark}\label{ModContD} (Mancini, 2005) Under \textbf{A2} we have a.s.
$$
\sup_{1\leq j \leq
n}\frac{|\Delta_{j}D^{(q)}|}{\sqrt{2hlog\frac{1}{h}}} \leq
K_{q}(\omega)<\infty,\quad q=1,2,
$$
where $K_q$ are finite random variables.
\end{Remark}

\begin{Theorem}\label{consistencyIA}(Estimation of the correlation between the continuous parts)
Let $(X^{(1)}_{t})_{t\in [0,T]}$ and $(X^{(2)}_{t})_{t\in [0,T]}$
two processes of the form (\ref{modello}). Assume
\textbf{A1}-\textbf{A4} and {\bf A6} are satisfied. Then
$$
\tilde{v}^{(n)}_{1,1}(X^{(1)},X^{(2)})_{T}\stackrel{P}{\longrightarrow}
\int_{0}^{T}\rho_{t} \sigma_{t}^{(1)}\sigma_{t}^{(2)}dt,
$$
as $n\rightarrow \infty$, where for $r$ and $l\in \N$
$$
\tilde{v}^{(n)}_{r,l}(X^{(1)},X^{(2)})_{T}
=h^{1-\frac{r+l}{2}}\sum_{j=1}^{n}(\Delta_{j}X^{(1)})^{r}1_{\{(\Delta_{j}X^{(1)})^{2}\leq
r_{h}\}}(\Delta_{j}X^{(2)})^{l}1_{\{(\Delta_{j}X^{(2)})^{2}\leq
r_{h}\}}.
$$\qed
\end{Theorem}


\n $v^{(n)}_{r,l}(X^{(1)},X^{(2)})_{T}
=h^{1-\frac{r+l}{2}}\sum_{j=1}^{n}(\Delta_{j}X^{(1)})^{r}(\Delta_{j}X^{(2)})^{l},
$ was used in \cite{BarNShe04} to estimate the covariation in the
case of diffusion processes. $\tilde
v^{(n)}_{r,l}(X^{(1)},X^{(2)})_{T}$ is a {\it threshold} modified
version
for the case of jump diffusion processes 
where 
we exclude from the sums the terms containing some jumps.

\begin{Remark}{\rm
An estimate of the sum of the co-jumps is obtained simply
subtracting the diffusion covariation estimator from the quadratic
covariation estimator. In fact
$$ \sum_{j=1}^{n}\Delta_{j}X^{(1)}\Delta_{j}X^{(2)}-
\tilde{v}^{(n)}_{1,1}(X^{(1)},X^{(2)})_{T}
\stackrel{P}\longrightarrow \sum_{0\leq s\leq T} \Delta
J_s^{(1)}\Delta J_s^{(2)},
$$
as $n\ri \infty$. Therefore an estimate of each $ \Delta
J_s^{(1)}\Delta J_s^{(2)}$ is obtained using
$$\Delta_{j}X^{(1)}\Delta_{j}X^{(2)}-\Delta_{j}X^{(1)}1_{\{(\Delta_{j}X^{(1)})^{2}\leq r_{h}\}}
\Delta_{j}X^{(2)}1_{\{(\Delta_{j}X^{(2)})^{2}\leq r_{h}\}},$$ with
$j$ such that $s\in ]t_{j-1}, t_j],$ whose limit for $h\ri 0$
coincides with the limit of
 $$\Delta_{j}X^{(1)}1_{\{(\Delta_{j}X^{(1)})^{2}> r_{h}\}}
\Delta_{j}X^{(2)}1_{\{(\Delta_{j}X^{(2)})^{2}> r_{h}\}}.
$$
 }\end{Remark}

\begin{Theorem}\label{CLTFiniteActivity} If $\tilde J_2^{(q)}\equiv 0,$ under the assumptions
\textbf{A1}-\textbf{A3}, and choosing $r_h$ as in \textbf{A6}, we
have
$$
{\cal
NB}(h):=\frac{\tilde{v}^{(n)}_{1,1}(X^{(1)},X^{(2)})_{T}-\int_{0}^{T}\rho_{t}
\sigma_{t}^{(1)}\sigma_{t}^{(2)}dt} {\sqrt
h\sqrt{\tilde{v}^{(n)}_{2,2}(X^{(1)},X^{(2)})_{T}-\tilde{w}^{(n)}(X^{(1)},X^{(2)})_{T}}}
\stackrel{d}{\longrightarrow} Z,
$$
where $Z$ has law ${\cal N}(0,1)$ and $$
\tilde{w}^{(n)}(X^{(1)}\!\!,X^{(2)})_{T}\!=\!h^{-1}\!\sum_{j=1}^{n-1}
\prod_{i=0}^{1}\Delta_{j+i}X^{(1)}1_{\{(\Delta_{j+i}X^{(1)})^{2}\leq
r_{h}\}}\prod_{i=0}^{1}\Delta_{j+i}X^{(2)}1_{\{(\Delta_{j+i}X^{(2)})^{2}\leq
r_{h}\}}.
$$

\end{Theorem}

\section{Main results}

In this paper we study the behavior of the normalized bias $ {\cal
NB}(h)$
 when infinite activity
jump components $\tilde J_2^{(q)}$ are included in the models
$X^{(q)}$. First we show that the standard error $$\sqrt
h\sqrt{\tilde{v}^{(n)}_{2,2}(X^{(1)},X^{(2)})_{T}-\tilde{w}^{(n)}(X^{(1)},X^{(2)})_{T}}$$
converges even in the present framework. We need the following
notations and remarks.

\begin{Remark}\label{4puntiPerDimoConsIA}[Remark 4.3 in \cite{GobMan06}] Under assumptions {\bf A2} and {\bf A5.2}
\begin{enumerate}
    \item If processes $a$ and $\sigma$ are càdlàg then, under {\bf A5}, a.s., for small $h$,
    $1_{\{(\Delta_{j}D^{(q)})^{2}>r_{h}\}}=0$, uniformly in $j$;
    \item Let us consider the
sequence $\tilde{v}_{1,1}^{(n)}, n\in\N$.  As long as $\tilde
J_2^{(q)}$ is a semimartingale, we can find a subsequence $n_k$
for which a.s., for large $k$, for all $j=1..n_k$, on $\{\Delta_j
X^{(q)}\leq 4 r(h_k)\}$ we have that
$$(\Delta \tilde J_{2,s})^2\leq 4r(h_k), \quad \forall s\in ]t_{j-1},
t_j].$$
    \item If $\tilde J_2^{(q)}$ is Lévy and independent of  $N^{(q)}$, and
    if $P\{ \Delta_j N \neq 0\}=O(h)$ as $h\ri 0$, then
     for any $j=1..n$, $nP\{ \DN\neq 0, \DXdq>\rh \}\ri 0 $ as $h\ri 0$.
\end{enumerate}
\end{Remark}


\n {\bf Notations}. For each $q=1,2$ we denote
$$\Delta_{j}\tilde J_{2 m}^{(q)} :=
\int_{t_{j-1}}^{t_{j}}\!\int_{|x|\leq 2\sqrt{r_{h}}}\ x\
\tilde{\mu}^{(q)}(dx,dt),\quad \Delta_{j}\tilde J_{2 c}^{(q)}
:=\int_{t_{j-1}}^{t_{j}}\int_{2\sqrt{r_{h}}<|x|\leq 1}x\
\nu^{(q)}(dx)dt$$ so that \beq\label{ScompDXddovepic}
\Delta_{j}\tilde{J}^{(q)}_{2}1_{\{|\Delta_{j}\tilde{J}_{2}^{(q)}|\leq
2\sqrt{r_{h}}\}}= \Delta_{j}\tilde J_{2 m}^{(q)} -
\Delta_{j}\tilde J_{2 c}^{(q)}. \eeq We also set
$$ \Delta_{j \star}H^{(q)}:= \Delta_{j}H^{(q)}1_{\{(\Delta_{j}X^{(q)})^{2}\leq r_{h}\}}$$
for any process $H^{(q)}$ (e.g. $H^{(q)}= Y^{(q)}$, or $H^{(q)}= \tilde J_2^{(q)}$ and so on).\\

Note that for each $q=1,2$
$$ E[ (\Delta_{j}\tilde J_{2 m}^{(q)})^2]=h\int_{|x|\leq 2\sqrt{r_h}} x^2 \nu^{(q)}(dx):=h\eta_q^2(2\sqrt{r_h}) \ri 0$$
as $h\ri 0,$ and under assumption {\bf A5} we have
\beq\label{DtildeJ2c}
  \Delta_{j}\tilde J_{2 c}^{(q)}= O\big(h(c - c r_h^{\frac{1-\alpha_q}{2}})\big).
  \eeq

%
%

\begin{Theorem}[standard error] Under the assumptions {\bf A1}-{\bf A6}, if
$\frac{hlog^{2}\frac{1}{h}}{r_{h}}\rightarrow 0$,
and $r_h=h^\beta$, $\beta\in ]0,1[$, then
$$
\tilde{v}^{(n)}_{2,2}(X^{(1)},X^{(2)})_{T}-\tilde{w}^{(n)}(X^{(1)},X^{(2)})_{T}\stackrel{P}{\longrightarrow}
\int_{0}^{T}(1+\rho_t ^{2})
(\sigma_{t}^{(1)})^{2}(\sigma_{t}^{(2)})^{2}dt
$$
as $n\rightarrow \infty$.
\end{Theorem}
\emph{\textbf{Proof}}. We prove that
$$
\tilde{v}^{(n)}_{2,2}(X^{(1)},X^{(2)})_{T}\stackrel{P}{\longrightarrow}
\int_{0}^{T}(2\rho_t^{2}+1)(\sigma_{t}^{(1)})^{2}(\sigma_{t}^{(2)})^{2}dt
$$
 and
$$
\tilde{w}^{(n)}(X^{(1)},X^{(2)})_{T}\stackrel{P}{\longrightarrow}
\int_{0}^{T}\rho_t^{2}(\sigma_{t}^{(1)})^{2}(\sigma_{t}^{(2)})^{2}dt.
$$
Note that a.s. for small $h$ that \beq\label{uno}
 1_{\{(\Delta_{j}X^{(q)})^{2}\leq r_{h}\}}=
 1_{\{(\Delta_{j}X^{(q)})^{2}\leq r_{h},(\Delta_{j}\tilde J_2^{(q)})^{2}\leq 4 r_{h} \}}+
 1_{\{(\Delta_{j}X^{(q)})^{2}\leq r_{h},(\Delta_{j}\tilde J_2^{(q)})^{2}> 4 r_{h} \}},\eeq
 and, trivially, we also have that
 \beq\label{due}
1_{\{|\Delta_{j}X^{(q)}|\leq\sqrt{r_{h}},|\Delta_{j}\tilde{J}_{2}^{(q)}|\leq
2\sqrt{r_{h}}\}}=
1_{\{|\Delta_{j}X^{(q)}|\leq\sqrt{r_{h}},|\Delta_{j}\tilde{J}_{2}^{(q)}|\leq
2\sqrt{r_{h}}, \DN^{(q)}=0\}}. \eeq

Let us now deal with $\tilde v_{2 2}^{(n)}$. As in the proof of
proposition 3.5 in \cite{GobMan06}  we can write
\beq\label{sommainiziale} \ba
\tilde{v}^{(n)}_{2,2}(X^{(1)},X^{(2)})_{T}-\int_{0}^
{T}(2\rho_t^{2}+1)(\sigma_{t}^{(1)})^{2}(\sigma_{t}^{(2)})^{2}dt=
\\
\\
\Big[h^{-1}\sum_{j=1}^{n}
(\Delta_{j\star}Y^{(1)})^{2}(\Delta_{j\star}Y^{(2)})^{2}
-\int_{0}^{T}(2\rho_t^{2}+1)(\sigma_{t}^{(1)})^{2}(\sigma_{t}^{(2)})^{2}dt\Big]+
\\
\\
h^{-1}\sum_{j=1}^{n}\big[(\Delta_{j\star}Y^{(1)})^{2}(\Delta_{j\star}\tilde{J}_{2}^{(2)})^{2}+
2(\Delta_{j\star}Y^{(1)})^{2}(\Delta_{j\star}Y^{(2)})(\Delta_{j}\tilde{J}_{2}^{(2)})+
(\Delta_{j\star}\tilde{J}_{2}^{(1)})^{2}(\Delta_{j\star}Y^{(2)})^{2}+
\\
\\
+(\Delta_{j\star}\tilde{J}_{2}^{(1)})^{2}(\Delta_{j\star}\tilde{J}_{2}^{(2)})^{2}+
2(\Delta_{j}\tilde{J}_{2}^{(1)})^{2}(\Delta_{j\star}Y^{(2)})(\Delta_{j\star}\tilde{J}_{2}^{(2)})+
2(\Delta_{j\star}Y^{(1)})(\Delta_{j\star}\tilde{J}_{2}^{(1)})(\Delta_{j}Y^{(2)})^{2}+
\\
\\
2(\Delta_{j\star}Y^{(1)})(\Delta_{j}\tilde{J}_{2}^{(1)})(\Delta_{j\star}\tilde{J}_{2}^{(2)})^{2}+
4(\Delta_{j\star}Y^{(1)})(\Delta_{j}\tilde{J}_{2}^{(1)})(\Delta_{j\star}Y^{(2)})(\Delta_{j}\tilde{J}_{2}^{(2)})
\big] := \sum\limits_{k=1}^9 I_k(h).\ea \eeq

\n The terms  of the right hand side within brackets are denoted
by $I_1(h)$ and can be split into two parts by adding and
subtracting the quantity $h^{-1}\sum_{j=1}^{n}
  (\Delta_{j}Y^{(1)})^{2}1_{\{(\Delta_{j}Y^{(1)})^{2}\leq
4r_{h}\}}(\Delta_{j}Y^{(2)})^{2}1_{\{(\Delta_{j}Y^{(2)})^{2}\leq
4r_{h}\}}$ in the following way
\begin{equation}\label{YY+resto}
\begin{array}{c}
|I_1(h)|=\big|h^{-1}\sum_{j=1}^{n}
  (\Delta_{j\star}Y^{(1)})^{2}(\Delta_{j\star}Y^{(2)})^{2}-
  \int_{0}^{T}(2\rho_t^{2}+1)(\sigma_{t}^{(1)})^{2}(\sigma_{t}^{(2)})^{2}dt\big|\leq
\\
\\
\!\!\!\!\!\big|h^{-1}\sum_{j=1}^{n}
  (\Delta_{j}Y^{(1)})^{2}1_{\{(\Delta_{j}Y^{(1)})^{2}\leq
4r_{h}\}}(\Delta_{j}Y^{(2)})^{2}1_{\{(\Delta_{j}Y^{(2)})^{2}\leq
4r_{h}\}}-\int_{0}^{T}(2\rho_t^{2}+1)(\sigma_{t}^{(1)})^{2}(\sigma_{t}^{(2)})^{2}dt\big|+
\\
\\
\!\!\!\!\big|h^{-1}\sum_{j=1}^{n}
  (\Delta_{j}Y^{(1)})^{2}(\Delta_{j}Y^{(2)})^{2}(1_{\{(\Delta_{j}X^{(1)})^{2}\leq
r_{h}\}}1_{\{(\Delta_{j}X^{(2)})^{2}\leq
r_{h}\}}-1_{\{(\Delta_{j}Y^{(1)})^{2}\leq
4r_{h}\}}1_{\{(\Delta_{j}Y^{(2)})^{2}\leq 4r_{h}\}})\big|
\end{array}
\end{equation}

\n The first term of the right hand side of (\ref{YY+resto}) tends
to zero in probability by proposition \ref{DenominPerCLTFA}.
Developing the second one we find that it is the sum of terms
which a.s. for small $h$ are zero because by remark
\ref{4puntiPerDimoConsIA} point 1 we have
\beq\label{IndDXpicmaDXdgr}1_{\{(\Delta_{j}X^{(q)})^{2}\leq
r_{h},(\Delta_{j}Y^{(q)})^{2}> 4r_{h}\}}\leq
1_{\{|\Delta_{j}\tilde{J}_{2}^{(q)}|> \sqrt{r_{h}}\}}\eeq and
\beq\label{IndDXgrmaDXdpic}1_{\{(\Delta_{j}X^{(q)})^{2}>
r_{h},(\Delta_{j}Y^{(q)})^{2}\leq 4r_{h}\}}\leq
1_{\{|\Delta_{j}D^{(q)}|> \frac{\sqrt{r_{h}}}{2}\}}+
1_{\{|\Delta_{j}\tilde{J}_{2}^{(q)}|>
\frac{\sqrt{r_{h}}}{2}\}}=1_{\{|\Delta_{j}\tilde{J}_{2}^{(q)}|>
\frac{\sqrt{r_{h}}}{2}\}},\eeq uniformly in $j$, so that the terms
containing $\Delta_jJ^{(q)}_{1}$ tends to zero by remark
\ref{4puntiPerDimoConsIA} point 3, whereas
$$
h^{-1}\sum_{j=1}^{n}
  (\Delta_{j}D^{(1)})^{2}(\Delta_{j}D^{(2)})^{2}1_{\{|\Delta_{j}\tilde{J}_{2}^{(1)}|>
\sqrt{r_{h}}\}}1_{\{|\Delta_{j}\tilde{J}_{2}^{(2)}|>
\sqrt{r_{h}}\}}\leq
$$
$$
K_{1}^{2}(\omega)K_{2}^{2}(\omega)hlog^{2}\frac{1}{h}\sum_{j=1}^{n}
  1_{\{|\Delta_{j}\tilde{J}_{2}^{(1)}|>
\sqrt{r_{h}}\}},
$$
which converges to zero in $L^1$
$$
E\big|hlog^{2}\frac{1}{h}\sum_{j=1}^{n}
  1_{\{|\Delta_{j}\tilde{J}_{2}^{(1)}|>
\sqrt{r_{h}}\}}\big|\leq nhlog^{2}\frac{1}{h}E\big[
  1_{\{|\Delta_{1}\tilde{J}_{2}^{(1)}|>
\sqrt{r_{h}}\}}\big]=
T\frac{hlog^{2}\frac{1}{h}}{r_{h}}\eta_{2}^{2}(1)\rightarrow 0
$$
The other terms in the right hand side of (\ref{sommainiziale})
tend to zero in probability. We only deal with $I_2, I_3, I_5,
I_8$ and $I_9$, the other ones being analogue. Note that for each
$q=1,2$
$$
E\left[\sup_{1\leq j\leq n}\frac{(\Delta_{j}\tilde{J}_{2}^{(q)})^2
1_{\{|\Delta_{j}\tilde{J}_{2}^{(q)}|\leq 2\sqrt{r_h}\}}}{h}\right]
\leq 2\sup_{1\leq j\leq
n}\frac{E(\Delta_{j}\tilde{J}^{(q)}_{2m})^{2}}{h}+ 2\sup_{1\leq
j\leq n}\frac{E(\Delta_{j}\tilde{J}^{(q)}_{2c})^{2}}{h}$$
$$=2\eta^{2}_q\Big(2\sqrt{r_h}\Big)+O\Big(2h(c-cr_{h}^{\frac{1-\alpha_{2}}{2}})^{2}\Big)=
O\big(h^{1+\beta(1-\alpha_q)}\big)$$ and $h^{1+\beta(1-\alpha_q)}$
tends to zero as $h\ri 0$. That is trivial if $\alpha_{q}\leq 1$;
however even if $\alpha_q$ belongs to $]1,2[$ it is ensured that
$1+\beta(1-\alpha_q)>0$, i.e. $\beta < \frac{1}{\alpha_q -1}$,
since $\beta<1$ while $\frac{1}{\alpha_q -1}>1.$ We have then that
as $h\ri 0$ \beq\label{quattro} \sup_{1\leq j\leq
n}\frac{(\Delta_{j}\tilde{J}_{2}^{(q)})^{2}
       1_{\{|\Delta_{j}\tilde{J}_{2}^{(q)}|\leq 2\sqrt{r_{h}}\}}}{h}\stackrel{P}{\longrightarrow} 0.\eeq

Now, by (\ref{uno}) and (\ref{due}) a.s. for small $h$
$$ |I_2+I_3+I_5+I_8 +I_9|\leq
h^{-1}\sum_{j=1}^{n}\Big|(\Delta_{j}D^{(1)})^{2}(\Delta_{j}\tilde{J}_{2}^{(2)})^{2}+
2(\Delta_{j}D^{(1)})^{2}(\Delta_{j}D^{(2)})(\Delta_{j}\tilde{J}_{2}^{(2)})
$$
$$
+(\Delta_{j}\tilde{J}_{2}^{(1)})^{2}(\Delta_{j}\tilde{J}_{2}^{(2)})^{2}
+2(\Delta_{j}D^{(1)})(\Delta_{j}\tilde{J}_{2}^{(1)})(\Delta_{j}\tilde{J}_{2}^{(2)})^2$$
$$+4(\Delta_{j}D^{(1)})(\Delta_{j}\tilde{J}_{2}^{(1)})(\Delta_{j}D^{(2)})(\Delta_{j}\tilde{J}_{2}^{(2)})
\Big| 1_{\{ (\DX^{(1)})^2\leq r_h\}} 1_{\{ (\DX^{(1)})^2\leq
r_h\}},
$$
and since $1_{\{ (\DX^{(q)})^2\leq r_h\}}=1_{\{ (\DX^{(1)})^2\leq
r_h,\Delta_{j}\tilde{J}_{2}^{(q)})^{2}\leq 2\sqrt{r_h}\}}+1_{\{
(\DX^{(1)})^2\leq r_h,\Delta_{j}\tilde{J}_{2}^{(q)})^{2}>
2\sqrt{r_h}\}}$ by (\ref{due}) the terms containing the indicator
of the set $\{\Delta_{j}\tilde{J}_{2}^{(q)})^{2}\leq
2\sqrt{r_h}\}$ are dominated by
$$ \sup_{1\leq j\leq n}
\frac{(\Delta_{j}\tilde{J}_{2}^{(2)})^{2}1_{\{|\Delta_{j}\tilde{J}_{2}^{(2)}|\leq
2\sqrt{r_{h}}\}}}{h}
 \Big[ \sum_{j=1}^n (\Delta_{j}D^{(1)})^{2}+
      \sum_{j=1}^n (\Delta_{j}\tilde{J}_{2}^{(1)})^2 +
 \sum_{j=1}^n(\Delta_{j}D^{(1)})(\Delta_{j}\tilde{J}_{2}^{(1)})\Big]
$$
$$
 +2\bar K^2 \mcq \sum_{j=1}^n(\Delta_{j}D^{(2)})(\Delta_{j}\tilde{J}_{2}^{(2)})
 $$
 $$
+4\sup_{1\leq j\leq n}
\frac{|\Delta_{j}\tilde{J}_{2}^{(1)}|1_{\{|\Delta_{j}\tilde{J}_{2}^{(1)}|\leq
2\sqrt{r_{h}}\}}}{\sqrt h}
 \sup_{1\leq j\leq n}
\frac{|\Delta_{j}\tilde{J}_{2}^{(2)}|1_{\{|\Delta_{j}\tilde{J}_{2}^{(2)}|\leq
2\sqrt{r_{h}}\}}}{\sqrt h} \sum_{j=1}^n
(\Delta_{j}D^{(1)})(\Delta_{j}D^{(2)}),
$$
where $\bar K:= \sqrt 2 (K_1\vee K_2)$. Each term tends to zero in
probability by (\ref{quattro}) and using that \beq\label{SomDDq}
\sum_{j=1}^n(\Delta_{j}D^{(q)})^{2} \stackrel{P}{\longrightarrow}
\int_0^T (\sigma_t^{(q)})^2 dt <\infty \mbox{ a.s. }, \eeq
$\sum_{j=1}^n(\Delta_{j}\tilde{J}_{2}^{(1)})^2\stackrel{P}{\longrightarrow}
T\int_{|x|\leq 1} x^2 \nu^{(1)}(dx) <\infty$ a.s.,
$\sum_{j=1}^n(\Delta_{j}D^{(q)})(\Delta_{j}\tilde{J}_{2}^{(q)})\stackrel{P}{\longrightarrow} [D^{(1)}, \tilde{J}_{2}^{(1)}]_T=0$ and\\
$\sum_{j=1}^n(\Delta_{j}D^{(1)})(\Delta_{j}D^{(2)})\stackrel{P}{\longrightarrow}
\int_0^T \rho_t\sigma_t^{(1)}\sigma^{(2)}_t dt< \infty $ a.s.,
where by $[M,N]$ we denote the quadratic covariation process
associated to two semimartingales $M$ and $N$ (see
\cite{ConTan04}).

It remains to consider
$$
h^{-1}\sum_{j=1}^{n}\Big|(\Delta_{j}D^{(1)})^{2}(\Delta_{j}\tilde{J}_{2}^{(2)})^{2}+
2(\Delta_{j}D^{(1)})^{2}(\Delta_{j}D^{(2)})(\Delta_{j}\tilde{J}_{2}^{(2)})
$$
$$
+(\Delta_{j}\tilde{J}_{2}^{(1)})^{2}(\Delta_{j}\tilde{J}_{2}^{(2)})^{2}
+2(\Delta_{j}D^{(1)})(\Delta_{j}\tilde{J}_{2}^{(1)})(\Delta_{j}\tilde{J}_{2}^{(2)})^2$$
$$
+4(\Delta_{j}D^{(1)})(\Delta_{j}\tilde{J}_{2}^{(1)})(\Delta_{j}D^{(2)})(\Delta_{j}\tilde{J}_{2}^{(2)})
\Big| 1_{\{ (\DX^{(1)})^2\leq
r_h,(\Delta_{j}\tilde{J}_{2}^{(1)})^{2}> 2\sqrt{r_h}\}} 1_{\{
(\DX^{(1)})^2\leq r_h,(\Delta_{j}\tilde{J}_{2}^{(2)})^{2}>
2\sqrt{r_h}\}}.
$$
Now observing that on $\{ (\DX^{(q)})^2\leq
r_h,(\Delta_{j}\tilde{J}_{2}^{(q)})^{2}> 2\sqrt{r_h}\}$, $q=1,2$,
we have $\{ (\Delta_j Y^{(q)})^2> r_h\}$, so that, a.s. for small
$h$
$$
1_{\{ (\DX^{(q)})^2\leq r_h,(\Delta_{j}\tilde{J}_{2}^{(q)})^{2}>
2\sqrt{r_h}\}}\leq 1_{\{ |\Delta_j J_{1}^{(q)}|>
\frac{\sqrt{r_h}}{2}\}}+1_{\{ |\Delta_j D^{(q)}|>
\frac{\sqrt{r_h}}{2}\}}\leq 1_{\{\Delta_j N^{(q)} \neq 0\}}
$$
by remark \ref{4puntiPerDimoConsIA} point 3 we note that all terms
tend to zero.

\n We can conclude that $I_2+I_3+I_5+I_8 +I_9
\stackrel{P}{\longrightarrow} 0$ as $h\ri 0$, and this concludes
the proof of the convergence of $\tilde v_{2 2}$.

Now, we show that
$\tilde{w}^{(n)}(X^{(1)},X^{(2)})_{T}\stackrel{P}{\longrightarrow}
\int_{0}^{T}\rho_t^{2}(\sigma_{t}^{(1)})^{2}(\sigma_{t}^{(2)})^{2}dt$.
Note that
$$
\left|h^{-1}\sum_{j=1}^{n-1}\prod_{q=1}^{2}\Delta_{j\star}X^{(q)}
\prod_{q=1}^{2}\Delta_{j+1,\star}X^{(q)}-\int_{0}^{T}\rho_t^{2}(\sigma_{t}^{(1)})^{2}(\sigma_{t}^{(2)})^{2}dt\right|
$$
is the sum of \beq\label{sette}
\left|h^{-1}\sum_{j=1}^{n-1}\prod_{q=1}^{2}\Delta_{j\star}Y^{(q)}
\prod_{q=1}^{2}\Delta_{j+1, \star}Y^{(q)}
-\int_{0}^{T}\rho_t^{2}(\sigma_{t}^{(1)})^{2}(\sigma_{t}^{(2)})^{2}dt\right|
\eeq and of other 15 terms of type $ h^{-1}\sum_{j=1}^n
\Delta_{j\star}M^{(1)}\Delta_{j+1,
\star}H^{(1)}\Delta_{j\star}M^{(2)}\Delta_{j+1,\star}H^{(2)}$
where (since $\DX^{(q)}= \Delta_j Y^{(q)}+ \DXd$ for each $q=1,2$)
both $M$ and $H$ can be $Y$ or $\tilde J_2$ and at least one
factor is the increment of one of the two $\tilde J_2^{(q)}$,
$q=1,2$. Each one of the 15 terms tends to zero in probability as
$h\ri 0$. In fact the terms where only one factor is the increment
of one of the $\tilde J_2^{(q)}$s are bounded by \beq\label{otto}
\sqrt{h^{-1}\sum_{j=1}^n (\Delta_{j+s}\tilde
J_2^{(q)})^21_{\{|\Delta_{j+s}\tilde J_2^{(q)}|\leq 2\sqrh\}}(
\Delta_{j+s}D^{(r)})^2} \sqrt{h^{-1}\sum_{j=1}^n (\Delta_{j+\bar
s}D^{(1)})^2(\Delta_{j+\bar s}D^{(2)})^2}, \eeq where $s=0$ or
$1$, $\bar s$ is 1 iff $s$ is 0 and $q, r\in \{1,2\}$. Using
(\ref{quattro}), (\ref{SomDDq}) and using that $h^{-1}\sum_{j=1}^n
(\Delta_{j+\bar s}D^{(1)})^2\cdot$ $(\Delta_{j+\bar s}D^{(2)})^2=
v_{22}(D^{(1)}, D^{(2)})_T$ converges to the a.s. finite
correlation term $\int_0^T
(1+2\rho^2_t)(\sigma_t^{(1)})^2(\sigma_t^{(2)})^2 dt$
(\cite{BarNShe04}, and cfr proposition \ref{DenominPerCLTFA}), we
reach that
(\ref{otto}) tends to zero in probability.\\
The terms containing two increments of kind $\tilde J_{j+s}^{(q)}$
are dominated in probability, thanks to (\ref{quattro}), by
$$ o(1)\sum_{j=1}^n \Delta_{j}D^{(r)}\Delta_{j+s}D^{(q)}\leq
o(1)\sqrt{\sum_{j=1}^n (\Delta_{j}D^{(r)})^2}\sqrt{\sum_{j=1}^n
(\Delta_{j+ s}D^{(q)})^2}\riP 0.$$ The terms containing three
increments of kind  $\tilde J_{j+s}^{(q)}$ are dominated by
$$ o(1)\sum_{j=1}^n \Delta_{j+u}\tilde J_2^{(r)}\Delta_{j+s}D^{(q)}\leq
o(1)\sqrt{\sum_{j=1}^n (\Delta_{j+u}\tilde
J_2^{(r)})^2}\sqrt{\sum_{j=1}^n (\Delta_{j+ s}D^{(q)})^2}\riP 0,$$
where $u, s\in \{0,1\}$. The unique term of type  (\ref{terms})
containing four increments of kind  $\tilde J_{j+s}^{(q)}$ is
simply dominated, thanks to (\ref{quattro}), by $ o(1)nh \ri 0$.

As for (\ref{sette}), adding and subtracting
$$
h^{-1}\sum_{j=1}^{n-1}\prod_{q=1}^{2}\Delta_{j}Y^{(q)}1_{\{(\Delta_{j}Y^{(q)})^{2}\leq
4r_{h}\}}
\prod_{q=1}^{2}\Delta_{j+1}Y^{(q)}1_{\{(\Delta_{j+1}Y^{(q)})^{2}\leq
4r_{h}\}},
$$
we obtain
$$
\Big|h^{-1}\sum_{j=1}^{n-1}\Big[\prod_{q=1}^{2}\Delta_{j}Y^{(q)}1_{\{(\Delta_{j}X^{(q)})^{2}\leq
r_{h}\}}
\prod_{q=1}^{2}\Delta_{j+1}Y^{(q)}1_{\{(\Delta_{j+1}X^{(q)})^{2}\leq
r_{h}\}}\Big]-\int_{0}^{T}\rho_t^{2}(\sigma_{t}^{(1)})^{2}(\sigma_{t}^{(2)})^{2}dt\Big|
$$
$$
\leq
\Big|h^{-1}\sum_{j=1}^{n-1}\Big[\prod_{q=1}^{2}\Delta_{j}Y^{(q)}1_{\{(\Delta_{j}Y^{(q)})^{2}\leq
4r_{h}\}}
\prod_{q=1}^{2}\Delta_{j+1}Y^{(q)}1_{\{(\Delta_{j+1}Y^{(q)})^{2}\leq
4r_{h}\}}\Big]-\int_{0}^{T}\rho_t^{2}(\sigma_{t}^{(1)})^{2}(\sigma_{t}^{(2)})^{2}dt\Big|
$$
$$
+
\Big|h^{-1}\sum_{j=1}^{n-1}\Delta_{j}Y^{(1)}\Delta_{j+1}Y^{(1)}\Delta_{j}Y^{(1)}
\Delta_{j}\tilde{J}^{(2)}_{2}\times
$$
$$
\times \big(1_{\{(\Delta_{j}X^{(1)})^{2}\leq
r_{h},(\Delta_{j+1}X^{(1)})^{2}\leq
r_{h},(\Delta_{j}X^{(2)})^{2}\leq
r_{h},(\Delta_{j+1}X^{(2)})^{2}\leq r_{h}\}}+
$$
$$
-1_{\{(\Delta_{j}Y^{(1)})^{2}\leq
4r_{h},(\Delta_{j+1}Y^{(1)})^{2}\leq
4r_{h},(\Delta_{j}Y^{(2)})^{2}\leq
4r_{h},(\Delta_{j+1}Y^{(2)})^{2}\leq 4r_{h}\}}\big)\Big|.
$$
The first term tends to zero in probability by theorem
\ref{DenominPerCLTFA}, whereas for the second one we note that
developing the difference of the two indicators we obtain a sum of
terms which are dominated by indicators as in
(\ref{IndDXpicmaDXdgr}) and (\ref{IndDXgrmaDXdpic}) and thus they
vanish a.s. for small $h$ (analogously as in
(\ref{YY+resto})).\qed

\vspace{1cm} Next we check the speed of convergence  to zero of
the estimation error
$\tilde{v}^{(n)}_{1,1}(X^{(1)},X^{(2)})_{T}-\int_{0}^{T}\rho_{t}
\sigma_{t}^{(1)}\sigma_{t}^{(2)}dt$. Within
$\tilde{v}^{(n)}_{1,1}(X^{(1)},X^{(2)})_{T}-\int_{0}^{T}\rho_{t}
\sigma_{t}^{(1)}\sigma_{t}^{(2)}dt$ it is the co-jumps term
$\sum_{s\leq t}\Delta \tilde J^{(1)}_{2 s}\Delta \tilde J^{(2)}_{2
s}$ to determine such a speed. However the speed of convergence of
such term depends both on the amount of jump activity of each
$\tilde J^{(q)}_{2}$ and on the dependence structure giving the
joint law $\big(\tilde J^{(1)}_{2}, \tilde J^{(2)}_{2}\big)(P)$.
We specialize our analysis to the case where $\tilde J^{(q)}_{2}$
have stable-like laws and the joint law is characterized by a
copula $C$ ranging in a given class.

\begin{figure}[h]\center{\begin{minipage}[c]{0.9\linewidth}
 \textbf{A7}
Assume $\alpha_q\in ]0,2[$ for each $q=1,2$. Consider (w.l.g.) $\alpha_1\leq \alpha_2$.\\
Each marginal law $\big(\tilde J^{(q)}_{2}\big)(P)$ has a
Stable-like density of the form $$\nu^{(q)}=c_q x^{-1-\alpha_q}
1_{\{x>0\}}+ d_q\  |x|^{-1-\alpha_q} 1_{\{x<0\}}.$$ For
simplicity, but w.l.g., we develop our proofs for the case where
each $\tilde J_2^{(q)}$ has only positive jump sizes, i.e.
$$\nu^{(q)}=c_q x^{-1-\alpha_q} 1_{\{x>0\}},$$
which have support $\R_+$.\\
We denote for each $q=1,2$ by \beq\label{tailq}
U_q(x):=\nu^{(q)}\big([x_q,+\infty[\big)= c_q
\frac{x_q^{-\alpha_q}}{\alpha_q}\eeq
 the tail integral of the marginal law of $\tilde J^{(q)}_{2}$.
 \end{minipage}}\end{figure}

 \begin{figure}[h]\center{\begin{minipage}[c]{0.9\linewidth}
 \textbf{A8}
The joint law $\big(\tilde J^{(1)}_{2}, \tilde
J^{(2)}_{2}\big)(P)$ has tail integrals given by
$$ U(x,y)= C_\gamma(U_1(x), U_2(y))$$
where $C_\gamma(u,v)$ is a Lévy copula (see \cite{ConTan04}) of
the form
 $$C_\gamma(u,v)=\gamma C_{\perp}(u,v) + (1-\gamma) C_{\parallel}(u,v),$$
where $C_{\perp}(u,v)= u 1_{\{v=\infty\}}+ v 1_{\{u=\infty\}}$ is
the independence copula, $C_{\parallel}(u,v)=u \wedge v$ is the
total dependence copula and $\gamma$ ranges in $[0,1]$.
\end{minipage}}\end{figure}

Such choices are quite representative since in fact many commonly
used models in finance (Variance Gamma model, CGMY model, NIG
model, etc.) have $\nu^{(q)}$ related to the ones in assumption
{\bf A7} in the sense that they are tempered stable processes
where the order of magnitude of the tail integrals as $x_q \ri 0$
is as for (\ref{tailq}). Moreover $C$ allows to range from a
framework of independent components to a framework where the
components are
completely positively monotonic.\\

\brem\label{RiduzCgenericaCparall} We need assumption {\bf A8} in
order to control the speed of convergence to zero of integrals
like $ \int_{0\leq x,y\leq \ep} xy d\nu(x,y) $, $ \int_{0\leq
x,y\leq \ep} x^2y^2 d\nu(x,y)$, where $\nu$ is the bivariate Lévy
measure of $(\tilde J_2^{(1)}, \tilde J_2^{(2)})$. Note that when
the copula within $\nu$ is the independence copula then both
integrals are zero so that under assumption {\bf A8}
$$ \int_{0\leq x,y\leq \ep} x^ky^k \nu(dx,dy) =
(1-\gamma)\int_{0\leq x,y\leq \ep} x^ky^k
dC_{\parallel}(U_1(x),U_2(y))$$ for $k=1,2$, and the speed is
given only by the complete dependence component. \erem

Now we compute the speed of convergence to zero of the small
co-increments of the two $\tilde J_2^{(q)}$.

\begin{Theorem}\label{CLTperCoIncrements}
Choose $r_h= h^\beta, \beta \in ]0,1[$ and $C_\gamma(u,v)\equiv
C_{\parallel}(u,v)$ (i.e. $\gamma=0$). Assume {\bf A1}-{\bf A8}.
Then
$$\frac{\sum_{j=1}^n \Delta_{j}\tilde{J}^{(1)}_{2}1_{\{(\Delta_{j}\tilde{J}_{2}^{(1)})^{2}\leq
4r_{h}\}}
  \Delta_{j}\tilde{J}^{(2)}_{2}1_{\{(\Delta_{j}\tilde{J}_{2}^{(2)})^{2}\leq
  4r_{h}\}}- nE[H'_{n1}]}{\sqrt{nVar(H'_{n1})}}\stackrel{d}{\longrightarrow} {\cal N}(0,1),$$ as $h\ri 0$, where
  for $j=1..n$
$$H_{nj}':=\Delta_{j}\tilde{J}^{(1)}_{2}1_{\{(\Delta_{j}\tilde{J}_{2}^{(1)})^{2}\leq 4r_{h}\}}
  \Delta_{j}\tilde{J}^{(2)}_{2}1_{\{(\Delta_{j}\tilde{J}_{2}^{(2)})^{2}\leq 4r_{h}\}}$$
is such that   as $h\ri 0$
$$E[H'_{nj}]=
O(h^{1+\beta\frac{\alpha_1+\alpha_2-\alpha_1\alpha_2}{2\alpha_1}})+h^2O\Big(
  (c-ch^{\beta\frac{1-\alpha_1}{2}})(c-ch^{\beta\frac{1-\alpha_2}{2}})\Big)$$
 and
  $$Var(H'_{nj})=
O(h^{2+\frac{\beta}{2}(4-\alpha_1-\alpha_2)})+O(
h^{1+\beta\frac{2\alpha_1+2\alpha_2-\alpha_1\alpha_2}{2\alpha_1}}).$$
\end{Theorem}



\dimo We use the Lindeberg-Feller theorem.  Using \textbf{A7} and
(\ref{ScompDXddovepic}) we have
  $$E[H_{nj}']=h\int_{]0,h^{\frac{\beta}{2}}]}\int_{]0,h^{\frac{\beta}{2}}]}xy\nu(dx,dy)+\Delta_{j}\tilde{J}^{(1)}_{2c}
  \Delta_{j}\tilde{J}^{(2)}_{2c}$$$$=h\int_{]0,h^{\frac{\beta}{2}}]}\int_{]0,h^{\frac{\beta}{2}}]}xydC_{\parallel}(U_1(x),U_2(y))+\Delta_{j}\tilde{J}^{(1)}_{2c}
  \Delta_{j}\tilde{J}^{(2)}_{2c}$$$$=h\int_{\frac{(h^{\frac{\beta}{2}})^{-\alpha_1}}{\alpha_1}\vee
  \frac{(h^{\frac{\beta}{2}})^{-\alpha_2}}{\alpha_2}}^{+\infty}U_{1}^{-1}(u)U_{2}^{-1}(u)du+\Delta_{j}\tilde{J}^{(1)}_{2c}
  \Delta_{j}\tilde{J}^{(2)}_{2c}$$
  $$=
  O(h^{1+\beta\frac{\alpha_1+\alpha_2-\alpha_1\alpha_2}{2\alpha_1}})+O(
  h(c-ch^{\beta\frac{1-\alpha_1}{2}}))O(
  h(c-ch^{\beta\frac{1-\alpha_2}{2}})).$$
  Moreover, since
  $$E[\Delta_{j}\tilde{J}^{(1)}_{2m}\Delta_{j}\tilde{J}^{(2)}_{2m}]^2=h\int_{]0,h^{\frac{\beta}{2}}]}
  \int_{]0,h^{\frac{\beta}{2}}]}x^2y^2\nu(dx,dy)$$
  $$+h^2\Big(\int_{]0,h^{\frac{\beta}{2}}]}\int_{]0,h^{\frac{\beta}{2}}]}x^2\nu(dx,dy)\Big)
  \Big(\int_{]0,h^{\frac{\beta}{2}}]}\int_{]0,h^{\frac{\beta}{2}}]}y^2\nu(dx,dy)\Big)$$
  $$+2h^2\Big(\int_{]0,h^{\frac{\beta}{2}}]}\int_{]0,h^{\frac{\beta}{2}}]}xy\nu(dx,dy)\Big)^2,$$
  $$E[(\Delta_{j}\tilde{J}^{(1)}_{2m})^2\ \Delta_{j}\tilde{J}^{(2)}_{2m}]=h\int_{]0,h^{\frac{\beta}{2}}]}
  \int_{]0,h^{\frac{\beta}{2}}]}x^2y\nu(dx,dy)$$
  and
  $$E[\Delta_{j}\tilde{J}^{(1)}_{2m}\ (\Delta_{j}\tilde{J}^{(2)}_{2m})^2]=h\int_{]0,h^{\frac{\beta}{2}}]}
  \int_{]0,h^{\frac{\beta}{2}}]}xy^2\nu(dx,dy),$$
we get
$$Var(H_{nj}')=
  O(h^{2+\frac{\beta}{2}(4-\alpha_1-\alpha_2)}) +O(+h^{1+\beta\frac{2\alpha_1+2\alpha_2-\alpha_1\alpha_2}{2\alpha_1}}).$$
Notice that $\forall \alpha_q\in ]0,2[$ we have
$\frac{\alpha_1+\alpha_2-\alpha_1\alpha_2}{2\alpha_1}>0$. Denote
$$H_{nj}=\frac{H_{nj}'-E[H_{nj}']}{\sqrt{nVar(H_{nj}')}}$$
the normalized versions of $H_{nj}'$. In order to verify the
Lindeberg condition we consider the following sets

$$\{|H_{nj}|>\eta\}=\left\{\frac{H_{nj}'-E[H_{nj}']}{\sqrt{n Var(H_{nj}')}}>\eta\right\}
=\Big\{|H_{nj}'-E[H_{nj}']|>\eta\sqrt{n Var(H_{nj}')}\Big\}.$$ We
show that in fact, for small $h$, $H_{nj}'\leq E[H_{nj}']+\sqrt{n
Var(H_{nj}')}\ \forall j$, thus $\{|H_{nj}|>\eta\}=\emptyset$.
Actually, after boring computations\footnote{These are available
if requested.} we reach that
$$E[H_{nj}']+\eta\sqrt{n\ Var(H_{nj}')}=O(h ^{1+\beta\frac{\alpha_1+\alpha_2-\alpha_1\alpha_2}{2\alpha_1}})+
O\Big(\sqrt{h^{1+\frac{\beta}{2}(4-\alpha_1-\alpha_2)}+
h^{\beta\frac{2\alpha_1+2\alpha_2-\alpha_1\alpha_2}{2\alpha_1}} }\
\Big)$$ as $h\ri 0$. Note that, using (\ref{ScompDXddovepic}) and
(\ref{DtildeJ2c}),
$$H_{nj}'= \Delta_{j}\tilde{J}^{(1)}_{2m}
\Delta_{j}\tilde{J}^{(2)}_{2m}-\Delta_{j}\tilde{J}^{(1)}_{2m}
\Delta_{j}\tilde{J}^{(2)}_{2c}-\Delta_{j}\tilde{J}^{(1)}_{2c}
\Delta_{j}\tilde{J}^{(2)}_{2m}+\Delta_{j}\tilde{J}^{(1)}_{2c}
\Delta_{j}\tilde{J}^{(2)}_{2c}$$$$=\Delta_{j}\tilde{J}^{(1)}_{2m}
\Delta_{j}\tilde{J}^{(2)}_{2m}-\Delta_{j}\tilde{J}^{(1)}_{2m}
O(h(c-ch^{\beta\frac{1-\alpha_2}{2}}))+$$$$-\Delta_{j}\tilde{J}^{(2)}_{2m}
O(h(c-ch^{\beta\frac{1-\alpha_1}{2}}))+O(
  h(c-ch^{\beta\frac{1-\alpha_1}{2}}))O(
  h(c-ch^{\beta\frac{1-\alpha_2}{2}})),$$
therefore $$
H_{nj}'=o\Big(E[H_{nj}']+\eta\sqrt{nVar(H_{nj}')}\Big)$$ as
$h\rightarrow 0$. Since $
  \frac{h^2(c-ch^{\beta\frac{1-\alpha_1}{2}})(
  c-ch^{\beta\frac{1-\alpha_2}{2}})}{
  \sqrt{h^{1+\frac{\beta}{2}(4-\alpha_1-\alpha_2)}}}\rightarrow 0$ it follows
  that
  $h^2(c-ch^{\beta\frac{1-\alpha_1}{2}})(
  c-ch^{\beta\frac{1-\alpha_2}{2}})=o\Big(E[H_{nj}']+\eta\sqrt{nVar(H_{nj}')}\Big)$. Moreover
for each $q=1,2$
$$\Delta_{j}\tilde{J}^{(q)}_{2m}
O(h(c-ch^{\beta\frac{1-\alpha_{q}}{2}}))$$$$=\Big(\Delta_{j}\tilde{J}^{(q)}_{2}1_{\{|\Delta_{j}\tilde{J}_{2}^{(q)}|\leq
2\sqrt{r_{h}}\}}+\int_{t_{j-1}}^{t_j}\int_{2\sqrt{r_h}\leq
|x|<1}x\nu^{(q)}(dx)dt\Big)O(h(c-ch^{\beta\frac{1-\alpha_{q}}{2}}))$$
$$\leq \big(2\sqrt{r_{h}}+O(h(c-ch^{\beta\frac{1-\alpha_{q}}{2}}))\big)O(h(c-ch^{\beta\frac{1-\alpha_{q}}{2}}))$$
$$=\big(O(h^{\frac{\beta}{2}})+O(h(c-ch^{\beta\frac{1-\alpha_{q}}{2}}))\big)O(h(c-ch^{\beta\frac{1-\alpha_{q}}{2}}))$$
$$=o(\sqrt{h^{1+\frac{\beta}{2}(4-\alpha_1-\alpha_2)}}),$$
so that as $h\rightarrow 0$ \beq\label{23}
\Delta_{j}\tilde{J}^{(q)}_{2m}
O(h(c-ch^{\beta\frac{1-\alpha_{q}}{2}}))=o\Big(E[H_{nj}']+\eta\sqrt{n
Var(H_{nj}')}\Big).\eeq Now using (\ref{ScompDXddovepic}) we can
write
$$\Delta_{j}\tilde{J}^{(1)}_{2m}\Delta_{j}\tilde{J}^{(2)}_{2m}\leq
    \Delta_{j}\tilde{J}^{(2)}_{2m}O(h^{\frac{\beta}{2}})+\Delta_{j}\tilde{J}^{(2)}_{2m}O(h(
  c-ch^{\beta\frac{1-\alpha_1}{2}})).$$
But
  $$E\left|\frac{(\Delta_{j}\tilde{J}^{(2)}_{2m})h^{\frac{\beta}{2}}}{\sqrt{h^{1+\frac{\beta}{2}(4-\alpha_1-\alpha_2)}}}\right|
  \leq \frac{h^{\frac{\beta}{2}}}{h^{\frac{1}{2}+\frac{\beta}{4}(4-\alpha_1-\alpha_2)}}
  \sqrt{E[(\Delta_{j}\tilde{J}^{(2)}_{2m})^2]}=
  \frac{h^{\beta-\frac{\beta\alpha_2}{4}}}{h^{\frac{\beta}{2}(4-\alpha_1-\alpha_2)}}\rightarrow
  0,$$ as $h \rightarrow 0$. It follows, using also (\ref{23}), that
  $\Delta_{j}\tilde{J}^{(1)}_{2m}\Delta_{j}\tilde{J}^{(2)}_{2m}=o\Big(E[H_{nj}']+\eta\sqrt{nVar(H_{nj}')}\Big)$.
Therefore for small $h$, uniformly on $j$, we have
$\{|H_{nj}|<\eta\}=\emptyset$ and the Lindeberg condition is
satisfied and the proof of theorem is complete.\qed

\section{Appendix}
\begin{Proposition}\label{DenominPerCLTFA}(Proposition 3.5 in \cite{GobMan06})
If $\tilde J_2^{(q)}\equiv 0,$ under the assumptions
\textbf{A1}-\textbf{A3}, and choosing $r_h$ as in \textbf{A5}, we
have
$$
\tilde{v}^{(n)}_{2,2}(X^{(1)},X^{(2)})_{T}\stackrel{P}{\longrightarrow}
\int_{0}^{T}(2\rho_{t}^{2}+1)(\sigma_{t}^{(1)})^{2}(\sigma_{t}^{(2)})^{2}dt,
$$
and
$$
\tilde{w}^{(n)}(X^{(1)},X^{(2)})_{T}\stackrel{P}{\longrightarrow}
\int_{0}^{T}\rho_{t}^{2}(\sigma_{t}^{(1)})^{2}(\sigma_{t}^{(2)})^{2}dt.
$$
\end{Proposition}

\begin{Theorem}[Lindeberg-Feller] Let $\{H_{nj},\ j=1,....,j_{n},\
n=1,2,....\}$ be a double array of r.v.s independent in each row
such that $EH_{nj}=0$ and $EH^{2}_{nj}=\sigma_{nj}^{2}<\infty$ for
each $n$ and $j$ and moreover
$\sum_{j=1}^{j_{n}}\sigma_{nj}^{2}=1$. Let $F_{nj}$ be the
distribution function of $H_{nj}$. In order that
\begin{enumerate}
    \item $\max_{1\leq j\leq j_{n}}P(|H_{nj}|>\epsilon)\rightarrow
    0$, $\forall \epsilon>0$,
    \item $\sum_{j=1}^{j_{n}}H_{nj}\stackrel{d}{\longrightarrow}
    {\cal N}(0,1)$,
\end{enumerate}
it is necessary and sufficient that for each $\eta>0$ that the
Lindeberg condition holds, i.e.
$$
\sum_{j=1}^{j_{n}}\int_{|x|>\eta}x^{2}F_{nj}(dx)=\sum_{j=1}^{j_{n}}EH_{nj}^{2}1_{\{|H_{nj}|>\eta\}}\rightarrow
0.
$$
\end{Theorem}

\end{document}